\documentclass[smallextended,referee,envcountsect,]{svjour3}
\smartqed
\usepackage{graphicx,xcolor}
\usepackage{amsmath,amssymb,mathtools,mathrsfs}
\usepackage{algorithm}
\usepackage{algorithmic}
\usepackage{booktabs}
\usepackage{hyperref}

\def\calM{\mathcal{M}}
\def\RX{\mathcal{R}}
\def\Rbkk{{\mathbb{R}^{k\times k}}}
\def\Rbnk{{\mathbb{R}^{n\times k}}}
\def\Rbnn{{\mathbb{R}^{n\times n}}}
\def\Rb{\mathbb{R}}

\def\grad{\mathrm{grad}}
\def\gradc{\mathrm{grad}_{gc}f}
\def\gradm{\mathrm{grad}_{\mathbf{M}_{X}}}

\def\St{\mathrm{St}}
\def\vecriz{\mathrm{vec}}
\def\iSt{\mathrm{iSt}}
\def\iStkn{\mathrm{iSt}(k,n)}

\def\sym{\mathrm{sym}}
\def\skew{\mathrm{skew}}

\def\tr{\mathrm{tr}}
\def\diag{\mathrm{diag}}
\def\SPD{\mathrm{SPD}}
\def\RZ{\mathcal{R}_X}

\def\exp{\mathrm{exp}}
\def\e{\mathrm{e}}
\def\d{\mathrm{d}} 
\def\ddt{\frac{\d}{\d t}}
\def\Df{\textsc{D}}
\def\Mx{\mathbf{M}_{X}}

\def\Mxa{\mathbf{M}_{X}^{-1}}
\def\gM{\mathrm{g}_{\mathbf{M}_{X}}}

\def\irm{\mathrm{i}}
\def\PZ{\mathscr{P}_{X}}

\newcommand{\skewset}{\mathrm{Skew}}
\newcommand{\symset}{\mathrm{Sym}}
\newcommand{\abs}[1]{\left|#1\right|}

\def\restrict#1{\raise-.5ex\hbox{\ensuremath|}_{#1}}

\begin{document}

\title{A generalized canonical metric for optimization on the indefinite Stiefel manifold}
\titlerunning{Generalized canonical metric for indefinite Stiefel manifold}


\author{Dinh Van Tiep \and Duong~Thi~Viet~An \and Nguyen~Thi~Ngoc~Oanh \and Nguyen~Thanh~Son}

\institute{Dinh Van Tiep \at 
	Faculty of Fundamental and Applied Sciences,\\
		Thai Nguyen University of Technology\\
		24131 Thai Nguyen, Vietnam\\ 
		tiep.dv@tnut.edu.vn
		\and
Duong Thi Viet An  \at
Department of Mathematics and Informatics,\\ 
Thai Nguyen University of Sciences\\
24118 Thai Nguyen, Vietnam\\
andtv@tnus.edu.vn
\and 
Nguyen Thi Ngoc Oanh \at
Department of Mathematics and Informatics,\\ 
Thai Nguyen University of Sciences\\
24118 Thai Nguyen, Vietnam\\
oanhntn@tnus.edu.vn
\and 
Nguyen Thanh Son, Corresponding author \at
Department of Mathematics and Informatics,\\ 
Thai Nguyen University of Sciences\\
24118 Thai Nguyen, Vietnam\\
ntson@tnus.edu.vn
}

\date{Received: date / Accepted: date}

\maketitle

\begin{abstract}
Various tasks in scientific computing can be modeled as an optimization problem on the indefinite Stiefel manifold.  
We address this using the Riemannian approach, which basically consists of equipping the feasible set 
with a Riemannian metric, preparing geometric tools such as orthogonal projections, formulae for Riemannian gradient, retraction and then extending an unconstrained optimization algorithm on the Euclidean space to the established manifold. The choice for the metric undoubtedly has a great influence on the method. In the previous work [D.V. Tiep and N.T. Son, A Riemannian gradient descent  method for optimization on the indefinite Stiefel manifold, arXiv:2410.22068v2[math.OC]], a tractable metric, which is indeed a family of Riemannian metrics defined by a symmetric positive-definite matrix depending on the contact point, has been used. In general, it requires solving a Lyapunov matrix equation every time
the gradient of the cost function is needed, which might significantly contribute to the computational cost. To address this issue, we
propose a new Riemannian metric for the indefinite Stiefel manifold. Furthermore, we construct the associated geometric structure, including a so-called quasi-geodesic and propose a retraction based on this curve. 
We then numerically verify the performance of the Riemannian gradient descent method associated with the new geometry and compare it with the previous work.
\end{abstract}


\section{Introduction}
Optimization of smooth functions of matrix variable $X$ subject to the constraint
of the form $\calM = \{X\ :\ X^TAX = J\}$, where $A$ is nonsingular and symmetric, and $J$ is symmetric,  involutory, i.e., satisfying $J^2 = I$, is an established line of research due to its wide range of applications. Various problems such as dimensionality reduction of data via principal component analysis, 
canonical correlation analysis, 
finding locality preserving projection, 
and computation of the generalized eigenvalues of a positive-definite matrix pencil can be formulated as an optimization problem with such a constraint; see, e.g., \cite{HardSST04,HeN03,Hote92,KovaV95,SameW82,TiepSon2025,YgerBGR12}. 

Additional conditions can be imposed on the matrices $A$ and $J$ upon the applications considered making the feasible set a differentiable manifold. When $A=I_n$ and $J=I_k, k\leq n$, the $n\times n$ and $k\times k$ identity matrices, respectively,  $\calM$ is termed the \emph{orthogonal Stiefel manifold} \cite{EdelAS98} 
$$\St(k,n) := \{X\in\Rbnk\ :\ X^TX = I_k\}.$$
Alternatively, if $A$ is an $n\times n$ real symmetric positive-definite (spd) matrix and $J=I_k$, the set is called the \emph{generalized Stiefel manifold} 
$$\St_A(k,n) := \{X\in\Rbnk\ :\ X^TAX = I_k\},$$
see, e.g., \cite{ShusA2023}. Although $\St(k,n)$ is a special case of $\St_A(k,n)$, one can see that $\St_A(k,n)$ can be turned into an orthogonal Stiefel manifold by the mapping $X\mapsto R^{-1}Y$ where $A=R^TR$ is 
the Cholesky factorization of $A$. 
A more general one in which $A$ is nonsingular, symmetric, and $J$ is symmetric, involutory, which covers the orthogonal and generalized Stiefel manifolds as special cases, is the \emph{indefinite Stiefel manifold}
\begin{equation}\label{eq:iSt}
\iSt_{A,J}(k,n) := \{X\in\Rbnk\ :\ X^TAX = J\}.
\end{equation}
This manifold was recently introduced in \cite{TiepSon2025}. Exploiting the rich geometry of the feasible set, a popular approach, known as \emph{Riemannian optimization}, is to extend well-known methods for unconstrained optimization in the Euclidean space to the case of these Riemannian manifolds once necessary geometric tools such as a Riemannian metric, an orthogonal projection, and a retraction are constructed. One is referred to classical works \cite{AdleDMMS02,EdelAS98}, monographs \cite{AbsiMS08,Boumal23,Sato2021,Udri1994}, and articles \cite{manopt,MishS2016,SatoA19,ShusA2023,TiepSon2025,WangDPY2024,YgerBGR12}, to name a few among many others dedicated to optimization methods on the orthogonal, generalized, and indefinite Stiefel manifolds.

Equipping a Riemannian metric for a differentiable manifold is an essential step in Riemannian optimization. Recall that a Riemannian metric is an inner product defined on each tangent space of the manifold depending smoothly on the contact point. It not only helps quantify objects such as the length of vectors, the measure of angles but also 
has a great effect on the behavior of numerical methods. Specifically, it is a 
means of defining the Riemannian gradient of a cost function, a crucial 
component in most of optimization methods. 

There are different ways to equip a manifold with a Riemannian metric depending on how the manifold is viewed. For some 
left- or right-invariant Riemannian metrics, and pseudo-Riemannian metrics equipped for the orthogonal Stiefel manifold as a 
quotient manifold of two orthogonal groups, one is referred to \cite[Subsect. 2.3]{EdelAS98}, \cite{JurdML2020,NishA2005};  
cf. \cite{BendZ21} for the symplectic Stiefel manifold. Here, we follow 
\cite{TiepSon2025} in which $\iSt_{A,J}(k,n)$ is shown to be an embedded submanifold of $\Rbnk$ of dimension \mbox{$nk-\frac{1}{2}k(k+1)$}. Then, a Riemannian metric is defined via an $n\times n$ spd $X$-smoothly depending matrix $\Mx$ in which the inner product of $Z_1,Z_2\in T_X\iSt_{A,J}(k,n)$ is determined as 
$$\gM(Z_1,Z_2):=\tr\left(Z_1^T\Mx Z_2\right).$$ 
This is not the most general way of defining a metric on a submanifold of 
$n\times k$ structured  matrices 
but it covers most cases in practice, e.g., the standard Euclidean metric and the canonical metric \cite{EdelAS98}. Using this formulation, the Riemannian gradient of smooth cost functions via orthogonally projecting onto the tangent space requires solving a Lyapunov matrix equation \cite[Prop. 3.3]{TiepSon2025}, which can be computationally expensive when $k$ is not small. It is 
surprising that in the case of the orthogonal Stiefel manifold, the Riemannian gradient with respect to the Euclidean and canonical metrics has 
closed-form expressions; see \cite[Sect. 2]{EdelAS98} and \cite[Subsect. 3.6]{AbsiMS08}. Moreover, in the generalized Stiefel manifold $\St_A(k,n)$ when $A$ is chosen to be the representing matrix $\Mx$, a closed-form Riemannian gradient is available \cite[Sect. 4]{WangDPY2024}. Meanwhile, for the symplectic Stiefel manifold, a relatively similar manifold to $\St_A(k,n),$ computing the Riemannian gradient needs solving a Lyapunov matrix equation under the Euclidean metric and a general tractable metric 
and does not when the canonical-like metric is used \cite{GSAS21a,GSAS21,GSS2024_2}. On the one hand, these facts raise a hope that one can design a Riemannian metric, most probably the canonical(-like) metric as in \cite{EdelAS98,GSAS21}, for the indefinite Stiefel manifold such that solving the Lyapunov matrix equation is not compulsory 
in order to compute the Riemannian gradient. 
Note, however, that 
the motivation for using the canonical(-like) metric in \cite{EdelAS98,GSAS21} is to balance the contribution of the two components 
of tangent vectors in the metric formulation
instead of obviating the necessity of solving the Lyapunov matrix equation as that for the 
metric proposed in this work. 
On the other hand, they inspire several 
questions:  (1) How are we on track 
to avoid solving the Lyapunov matrix equations to compute the Riemannian gradient? (2) How to explain the different situations for orthogonal, generalized, and symplectic Stiefel manifolds mentioned above? And (3) how does an optimization method such as Riemannian gradient descent behave using the geometric quantities under the proposed metric?

\paragraph{Contribution} It is the main motivation source of the present paper to answer these questions. Namely, by converting the Lyapunov matrix equation into a linear equation using the Kronecker product and vectorization operator and by asking that the resulting equation is solving-free, we first achieve a general one-parameter family of metrics that fulfill our expectation; the canonical(-like) metric in \cite{EdelAS98,GSAS21} certainly belongs to this class. Second, we explain in detail why in the case of the orthogonal Stiefel manifold, the Lyapunov matrix equation does not appear 
during the computation of the Riemannian gradient by examining the formulation of the orthogonal projection onto the tangent space under the tractable metric. Third, we go further by characterizing the normal space to the manifold, deriving 
closed-form expressions for the orthogonal projections onto the tangent and normal spaces, and the Riemannian gradient of the cost function. Moreover, we construct a so-called quasi-geodesic curve associated with the obtained metric which in turn allows us to propose a retraction. 
Finally, these tools are invoked in a Riemannian gradient descent 
scheme for running several numerical examples.

\paragraph{Outline}  After recalling some main facts about the indefinite Stiefel manifold in Section~\ref{sec:Geometry}, we derive the generalized canonical metric that yields 
closed-form expressions for the orthogonal projections and the Riemannian gradient in Section~\ref{sec:metric}. An explanation for the aforementioned situations in orthogonal, generalized, and symplectic Stiefel manifolds is also given. In Section~\ref{sec:quasigeo}, we construct a quasi-geodesic 
which is used to propose a retraction in Section~\ref{sec:rectraction}. 
In Section~\ref{sec:NumerExam}, using a Riemannian gradient descent algorithm, we verify the theoretical findings via several examples. Finally, a conclusion is given in Section~\ref{sec:Concl}.

\paragraph{Notation} 
We finish this section by introducing the notations. We use $\Omega^T, \tr(\Omega),$ $\|\Omega\|_F,$ and $\|\Omega\|$ 
to signify the transpose, trace, Frobenius norm, and spectral norm of a matrix $\Omega,$ respectively. When $\Omega$ is a square matrix, 
its skew-symmetric and symmetric parts are denoted by 
$\skew(\Omega) \coloneqq \frac{1}{2}(\Omega-\Omega^{T})$ and $\sym(\Omega) \coloneqq \frac{1}{2}(\Omega + \Omega^{T}),$ respectively. The set of skew-symmetric, the set of symmetric, and the set of symmetric positive-definite matrices of size $p\times p$ are respectively denoted by $\skewset(p)$, $\symset(p)$, and 
$\SPD(p)$. We write $\diag(\cdot,\ldots,\cdot)$ to mean the (block) diagonal matrix. The matrix exponential is represented by either $\e^{\cdot}$ or $\exp(\cdot).$ In a linear space, $\langle\cdot,\cdot\rangle$ denotes the standard Euclidean inner product; specifically, $\langle P,Q\rangle = \tr(P^TQ)$ for $P, Q\in{\Rbnk}$. The first-order derivative with respect to $t$ of a smooth mapping $G(t)$ is represented as $\ddt G$ or $\dot{G}(t),$ and $\ddot{G}(t)$ means its second-order derivative. For a map $F:\mathcal{V}_1\longrightarrow\mathcal{V}_2$ between two normed vector spaces $\mathcal{V}_1$ and $\mathcal{V}_2,$ let $o(\cdot)$ denote the little-$o$ notation, the Fréchet derivative of $F$ at $X\in\mathcal{V}_1$ is a linear operator $\Df F(X):\mathcal{V}_1\longrightarrow\mathcal{V}_2$ fulfilling that $F(X+Z) = F(X)+\Df F(X)[Z]+o(\| Z\|),$ for any $Z\in\mathcal{V}_1$. 
\section{Preliminaries}
\label{sec:Geometry}
In this section, we recall basic facts about the indefinite Stiefel manifold $\iSt_{A,J}(k,n)$ which can also be shortened as $\iSt_{A,J}$ or $\iSt(k,n)$ if the dropped information is either inessential or clear from the context. For their derivations, the reader is referred to \cite{TiepSon2025}. 

Given matrix $A$ as introduced in the beginning, let $\irm_+(A)$ and $\irm_-(A)$ respectively denote the number of positive and negative eigenvalues of $A$. It was stated in \cite[Lem. 2.1]{TiepSon2025}  that a necessary and sufficient condition for the nonemptiness of the feasible set $\iSt_{A,J}$ is that 
\begin{equation}\label{eq:lem_nonempty}
\irm_+(J) \leq \irm_+(A) \mbox{ and }  \irm_-(J) \leq \irm_-(A).
\end{equation} 
We always assume that \eqref{eq:lem_nonempty} holds when working with this manifold. 

Then, it was shown in \cite[Prop. 2.2]{TiepSon2025} that the set 
$\iSt_{A,J}$ defined in \eqref{eq:iSt} is a closed, embedded submanifold of $\Rbnk$ of dimension $nk-\frac{1}{2}k(k+1)$. It was also noted that the indefinite Stiefel manifold is generally unbounded.

For 
$X\in\iSt_{A,J}(k,n)$, let us denote by $X_{\perp}\in\mathbb{R}^{n\times(n-k)}$ a full-rank matrix satisfying 
$X^{T}X_{\perp}=0\in\mathbb{R}^{k\times(n-k)}$ 
and set $E\coloneqq [X \;\; A^{-1}X_{\perp}]\in \Rbnn.$ The following statements hold due to \cite[Lem. 2.3]{TiepSon2025}:
\begin{itemize}
\item[i)] $E$ is nonsingular and
\begin{equation}\label{eq:Einverse}
E^{-1}=
\begin{bmatrix}
JX^{T}A \\ 
(X_{\perp}^{T}A^{-1}X_{\perp})^{-1}X_{\perp}^{T} 
\end{bmatrix}.
\end{equation}
\item[ii)] For any $Z\in\Rbnk,$ there are $W\in\Rbkk$ and $K\in\mathbb{R}^{(n-k)\times k} $ such that
\begin{equation}\label{eq:Zdecompose}
Z = XW + A^{-1}X_{\perp}K,
\end{equation}
where $W = JX^{T}AZ, K = (X_{\perp}^{T}A^{-1}X_{\perp})^{-1}X_{\perp}^{T}Z$. It is immediately implied from \eqref{eq:Zdecompose} that $A^{-1}X_{\perp}K = Z- XW$ is independent of $X_{\perp}.$
\item[iii)] There holds
\begin{equation}\label{eq:In_1}
I_{n} = XJX^TA + A^{-1}X_{\perp}(X_{\perp}^{\top}A^{-1}X_{\perp})^{-1}X_{\perp}^T.
\end{equation}
\end{itemize} 

%
%
Next, recall that the tangent space at $X\in  \iStkn$ is the set
\begin{equation*}
T_{X}\iStkn = \left\lbrace \dot{\varphi}(0) \;  :\;  \varphi: (a,b)\ni 0 \longrightarrow \iStkn \text{ is smooth and }\varphi(0)=X \right\rbrace.
\end{equation*} 
In \cite[Prop. 2.4]{TiepSon2025}, it has been shown that
\begin{align} \notag 
T_{X}\iStkn & = \left\{Z\in\Rbnk: Z^{T}AX + X^{T}AZ = 0\right\} \\
& = \left\{XW + A^{-1}X_{\perp}K: JW\in\skewset(k),    K\in\mathbb{R}^{(n-k)\times k} \right\}. \label{eq:eq31b}
\end{align}


In \cite[Sect. 3]{TiepSon2025}, a so-called tractable metric $\gM$ defined via an $X$-smoothly dependent matrix $\Mx\in\SPD(n)$ was proposed. Namely, for any pair $Z_1,Z_2 \in T_X\iStkn,$ $ \langle Z_1,Z_2\rangle_{\gM}\coloneq \tr\left(Z_1^T\Mx Z_2\right)$. 
With this metric, for any $X\in\iSt(k,n)$, the ambient space $\Rbnk$ is decomposed as the direct sum of the tangent space and the normal space, 
$$\Rbnk = T_X\iStkn \overset{\perp_{\gM}}{\bigoplus}T_X\iStkn^{\perp}_{\gM}.$$
This feature motivates 
definitions of the orthogonal projections onto the tangent space and onto the normal space which 
associate each $Y\in\Rbnk$ with a unique $\PZ(Y)\in T_X\iStkn$ and a unique $\PZ^{\perp}(Y)\in T_X\iStkn^{\perp}_{\gM}$, respectively, such that
\begin{equation}\label{eq:orthogonalsum}
Y = \PZ(Y) + \PZ^{\perp}(Y).
\end{equation}
The computation of $\PZ(Y)$ is essential in formulating the Riemannian gradient. According to \cite[Prop. 3.2]{TiepSon2025}, $\PZ(Y) = Y-\Mx^{-1}AXU_{X,Y}$ where $U_{X,Y}$ 
is the solution to the Lyapunov matrix equation

\begin{equation}\label{eq:Lya_Uxy}
2\sym(X^TAY) = X^TA\Mxa A X U_{X,Y}+U_{X,Y}X^TA\Mxa A X.
\end{equation}
Based on this formula, the Riemannian gradient was given in \cite[Prop. 3.3]{TiepSon2025} by 
\begin{equation}\label{eq:gradMx_formula}
\gradm f(X) = \Mxa\nabla\bar{f}(X) - \Mxa AX U_{\bar{f}},
\end{equation}
where $\bar{f}$ is some smooth extension of $f$ on the ambient space $\Rbnk$, $U_{\bar{f}}$ 
is the unique solution to the Lyapunov matrix equation with unknown $U$: 
\begin{equation}\label{eq:Lyapunov_gradMx}
(X^{T}A\Mxa A X)U+U(X^{T}A\Mxa A X) = 2\sym\left[ X^{T}A\Mxa\nabla\bar{f}(X)\right].
\end{equation}
Ignoring the computational cost for $\Mxa(\cdot)$ by, e.g., considering the Euclidean metric, one can show that the total cost for computing the Riemannian gradient $\gradm f(X)$ in \eqref{eq:gradMx_formula} is about $2n^2k+6nk^2-nk-2k^2 + 20k^3$ flops. The last term in this expression is contributed by solving the equation \eqref{eq:Lyapunov_gradMx} with an spd dense matrix coefficient using the Bartels-Stewart algorithm \cite{BartS72}. When $k$ is large, for instance, the same as $n$, it dominates the total cost. Therefore, avoiding solving this equation in such a case most probably helps accelerate the optimization process.

\section{A generalized canonical metric on $\iStkn$}\label{sec:metric} 
As introduced, we will determine a form of $\Mx$ to eliminate the need of 
solving the Lyapunov 
equation~\eqref{eq:Lyapunov_gradMx}.  
To this end, we simplify the coefficient $X^TA\Mxa A X$ of the equation. The  identity $E^TAX = \left[\begin{smallmatrix}
J\\
0
\end{smallmatrix}\right]$ suggests 
considering $\Mxa\in \SPD(n)$ as 
\begin{equation*}
\Mxa = E\begin{bmatrix}
\Gamma_1&\Gamma_2\\
\Gamma_2^T&\Gamma_3
\end{bmatrix}E^T, \ \Gamma_1\in\SPD(k),\;  \Gamma_2\in\Rb^{k\times(n-k)},\; \Gamma_3\in\SPD(n-k). 
\end{equation*}
This results in $X^TA\Mxa A X=J\Gamma_1 J.$ Consequently, 
the equation \eqref{eq:Lya_Uxy} reduces to
\begin{equation}\label{eq:Lya_Uxy_2}
2\sym(X^TAY) = (J\Gamma_1 J)U_{X,Y}+U_{X,Y}(J\Gamma_1 J).
\end{equation}
For greater detail, let us convert this matrix equation into a linear equation using vectorization operator and Kronecker product \cite{HornJ91}. Namely, for a matrix $M=(m_{ij})\in \Rbnk$, $\vecriz(M) \in \mathbb{R}^{nk}$ is generated by stacking the columns of $M$ to make a column vector, and the Kronecker product of $M$ 
and $N\in \mathbb{R}^{p\times q}$, denoted by $M\otimes N\in\mathbb{R}^{np\times kq}$, is the matrix consisting of $nk$ block matrices $m_{ij}N$ for $i=1,\ldots,n, j = 1,\ldots,k$. In view of \cite[Lem. 4.3.1]{HornJ91}, 
\eqref{eq:Lya_Uxy_2} is equivalent to
\begin{equation*}
(I_k\otimes (J\Gamma_1 J) + (J\Gamma_1 J)\otimes I_k)\vecriz(U_{X,Y}) = \vecriz(2\sym(X^TAY)).
\end{equation*}
We say that a square system of linear equations is \emph{solving-free} if its coefficient is a scalar matrix, i.e., the product of a (nonzero) scalar and the identity matrix. Thus, its solution is merely the product of the reciprocal of the scalar and the right-hand side. Accordingly, a sufficient condition for 
obviating the need to solve 
the Lyapunov matrix equation \eqref{eq:Lya_Uxy} 
is that $I_k\otimes (J\Gamma_1 J) + (J\Gamma_1 J)\otimes I_k$ 
has 
the form $2\rho I_{k^2}$, where $\rho$ is a nonzero constant parameter. It follows that $J\Gamma_1 J = \rho I_{k}$ which, together with the condition $J^2=I_k$, yields $\Gamma_1 = \rho I_{k}$. For the symmetric positive-definiteness of $\Mxa$, we set $\rho >0$. One can verify that in this case, the solution to \eqref{eq:Lya_Uxy_2} is 
$U_{X,Y} = \frac{1}{\rho}\sym(X^TAY).$
Furthermore, we set $\Gamma_2 = 0$ for the ease of argument. Consequently, $\Mxa$ takes the form
\begin{equation} \label{eq:explicitMx_2}
\Mxa = E\begin{bmatrix}
\rho I_k&0\\
0&\Gamma_3
\end{bmatrix}E^T = \rho XX^T +A^{-1}X_\perp\Gamma_3X_\perp^TA^{-1}.
\end{equation}
Moreover, taking the formula \eqref{eq:Einverse} for the inverse of $E$ into account, we also have
\begin{align}\notag
\Mx &=\begin{bmatrix}
JX^{T}A \\ 
(X_{\perp}^{T}A^{-1}X_{\perp})^{-1}X_{\perp}^{T} 
\end{bmatrix}^{T}
\begin{bmatrix}
\frac{1}{\rho}I_k& 0\\
0&\Gamma_3^{-1}\end{bmatrix}
\begin{bmatrix}
JX^{T}A \\ 
(X_{\perp}^{T}A^{-1}X_{\perp})^{-1}X_{\perp}^{T}\end{bmatrix} \\
&=\dfrac{1}{\rho}AXX^TA + X_{\perp}\left(X_{\perp}^TA^{-1}X_{\perp}\right)^{-1}\Gamma_3^{-1}\left(X_{\perp}^TA^{-1}X_{\perp}\right)^{-1}X_{\perp}^T. 
\label{eq:M_X_general}
\end{align}

In the framework of the tractable metric whose representing matrix has just been obtained, the targeted \textit{generalized canonical metric} applied to two tangent vectors $Z_i = XW_i + A^{-1}X_{\perp}K_i$, 
$i=1,2,$ at any $X\in\iStkn$ is 
\begin{equation}\label{eq:canon_metric}
g_{\rho,X_{\perp}}(Z_1,Z_2) \coloneqq \tr\left(Z_1^T\Mx Z_2\right) = \dfrac{1}{\rho}\tr(W_1^TW_2) + \tr(K_1^T\Gamma_3^{-1}K_2).
\end{equation}
Obviously, $g_{\rho,X_{\perp}}$ is a non-degenerate bilinear form defined on 
$T_{X}\iStkn$ and dependent upon 
both $\rho$ and $X_{\perp}$. To make this a Riemannian metric, it is essential to require that $\Mx$ depends only on $X\in\iStkn$ and this dependence is smooth. Apparently, we can consider $\Gamma_3$ as an independent parameter like $\rho$ and smoothly represent $X_\perp$ in terms of $X$ (or do this for both $\Gamma_3$ and $X_\perp$). Nevertheless, we do not know an efficient and explicit formula for this representation, if it exists.  Following \cite{EdelAS98,GSAS21}, we impose a relation on these two components and obtain the desired metric. Namely, we can choose 
\begin{equation}\label{eq:Xperp_1}
\left(X_\perp^TX_\perp\right)^{-1} = \Gamma_3
\end{equation}
then the second term of \eqref{eq:M_X_general} becomes
\begin{equation*}
X_{\perp}\left(X_{\perp}^TA^{-1}X_{\perp}\right)^{-1}\left(X_\perp^TX_\perp\right)\left(X_{\perp}^TA^{-1}X_{\perp}\right)^{-1}X_{\perp}^T
=\left(A(I_n-XJX^TA)\right)^2
\end{equation*}
thanks to the identity \eqref{eq:In_1}. Thus, the representing matrix is
\begin{equation}\label{eq:Mx_choice1}
\Mx = \frac{1}{\rho}AXX^TA + \left(A\left(I_n-XJX^TA\right)\right)^2.
\end{equation}
Alternatively, we can take
\begin{equation}\label{eq:Xperp_2}
\left(X_{\perp}^TA^{-1}X_{\perp}\right)^{-1}X_\perp^TX_{\perp}\left(X_{\perp}^TA^{-1}X_{\perp}\right)^{-1} = \Gamma_3.
\end{equation}
Accordingly, it follows that
\begin{align*}
X_{\perp}\left(X_{\perp}^TA^{-1}X_{\perp}\right)^{-1}\Gamma_3\left(X_{\perp}^TA^{-1}X_{\perp}\right)^{-1}X_{\perp}^T &= X_{\perp}\left(X_{\perp}^TX_{\perp}\right)^{-1}X_{\perp}^T\\
&= I_n - X\left(X^TX\right)^{-1}X^T.
\end{align*} 
Thus, 
\begin{equation*}
\Mx = \frac{1}{\rho}AXX^TA + I_n - X\left(X^TX\right)^{-1}X^T.
\end{equation*}
\begin{remark}
\begin{itemize}
	\item The relations \eqref{eq:Xperp_1} and \eqref{eq:Xperp_2} are only the means to obtain $\Mx$. As we will see later, we do not have to construct $\Gamma_3$ or $X_\perp$.
\item One can see that when $\rho = 2$ and $\Gamma_3 = I_{n-k}$, the formulation \eqref{eq:Mx_choice1} reduces to the canonical metric given in \cite[Eq.~(2.39)]{EdelAS98}. Hence, it is named the generalized canonical metric. Moreover, we use this value for $\rho$ in numerical tests unless stated otherwise.	
\item The fact that the generalized canonical metric \eqref{eq:canon_metric} is independent of $X_\perp$ upon the choices for $\Gamma_3$ can also be derived using \eqref{eq:Zdecompose}.  
Indeed, the dependence of the metric on $X_\perp$ is only on the second term which is now eliminated by using either \eqref{eq:Xperp_1} or \eqref{eq:Xperp_2}. 
Namely, for the first choice, $$K^T_1\Gamma_3^{-1}K_2 = K^T_1X_\perp^T X_\perp K_2 = \left(X_\perp K_1\right)^TX_\perp K_2$$ 
is independent of $X_\perp$ because $A^{-1}X_\perp K_j$ is so for $j=1,2,$ according to \eqref{eq:Zdecompose}. For the second one, we obtain that
\begin{align*}
K^T_1\Gamma_3^{-1}K_2 &= K^T_1\left(X_{\perp}^TA^{-1}X_{\perp}\right)\left(X_\perp^TX_{\perp}\right)^{-1}\left(X_{\perp}^TA^{-1}X_{\perp}\right)K_2\\
&=Z_1^TX_\perp\left(X_{\perp}^TA^{-1}X_{\perp}\right)^{-1}\left(X_{\perp}^TA^{-1}X_{\perp}\right)\left(X_\perp^TX_{\perp}\right)^{-1}\left(X_{\perp}^TA^{-1}X_{\perp}\right)\\
&\qquad\qquad\qquad\times\left(X_{\perp}^TA^{-1}X_{\perp}\right)^{-1}X_\perp^T Z_2\\
&= Z_1^TX_\perp\left(X_\perp^TX_{\perp}\right)^{-1}X_\perp^TZ_2\\
&= Z_1^T\left(I_n -X\left(X^TX\right)^{-1}X^T\right) Z_2.
\end{align*}
\end{itemize}
\end{remark}

\subsection{Geometry under the generalized canonical metric}
For any $X\in\iStkn,$ the matrix $\Mx$ defined by \eqref{eq:M_X_general} is symmetric positive-definite, so we can extend the inner product $g_{\rho,X_{\perp}}$ on $T_X\iStkn$ to an inner product on $\Rbnk$ which is given by $g_{\rho,X_{\perp}}(Z_1,Z_2)=\tr(Z_1^T\Mx Z_2),$ for all $Z_1,Z_2\in\Rbnk.$ For the moment, neither \eqref{eq:Xperp_1} nor \eqref{eq:Xperp_2} is imposed for the sake of generality. 
The normal space to $\iStkn$ at $X$ with respect to $g_{\rho,X_{\perp}}$ is defined by 
$$T_X\iStkn_{\rho}^{\perp}\coloneqq \left\lbrace N\in\Rbnk: g_{\rho,X_{\perp}}(N,Z) = 0 \text{ for all } Z\in T_X\iStkn \right\rbrace.$$ 
The following statement characterizes the elements of this space. Although the metric \eqref{eq:M_X_general} depends on $X_\perp$, the elements of the normal space, as we will see, do not. 
\begin{proposition}\label{prop:normsp2}
$T_X\iStkn_{\rho}^{\perp} = \left\lbrace XU: JU\in\symset(k) \right\rbrace.$ 
\end{proposition}

{\it Proof}
Let $N = XU_N + A^{-1}X_{\perp}K_N$ 
with $U_N\in\Rbkk,$ $K_N\in\mathbb{R}^{(n-k)\times k}$ as in \eqref{eq:Zdecompose}. 
For each $Z\in T_X\iStkn,$ we have $Z = XW_Z + A^{-1}X_{\perp}K_Z,$ with $JW_Z\in\skewset(k),K_Z\in\mathbb{R}^{(n-k)\times k}$ as in \eqref{eq:eq31b}. By \eqref{eq:canon_metric}, we have
\begin{align*}
g_{\rho,X_{\perp}}(N,Z) &= \dfrac{1}{\rho}\tr(U_{N}^TW_Z) + \tr(K_N^T\Gamma_3^{-1} K_Z)\\
&=\dfrac{1}{\rho}\tr\left[(JU_{N})^T(JW_Z)\right] + \tr(K_N^T\Gamma^{-1}_3 K_Z).
\end{align*}
Therefore, $g_{\rho,X_{\perp}}(N,Z)=0$, for any $JW_Z\in\skewset(k)$ and any $K_Z\in\mathbb{R}^{(n-k)\times k},$ 
if and only if
$K_N = 0$ and $JU_N\in\symset(k)$ since $\skewset(k)$ is the orthogonal complement of $\symset(k)$ with respect to the Euclidean product on $\Rbnk$. 
Therefore, $N = XU_N$ with $JU_N\in\symset(k)$ if and only if $N\in T_X\iStkn_{\rho}^{\perp}$ which is also the expected characterization. 
\qed

Next, we present formulae for the orthogonal projections $\PZ$ and 
$\PZ^{\perp}$ with respect to the inner product $g_{\rho,X_{\perp}}$ onto $T_X\iStkn$ and $T_X\iStkn^{\perp}_{\rho},$ respectively. 

\begin{proposition}\label{prop:projectiongrho}
Given any $Y\in\Rbnk,$ the orthogonal projections of  $Y$ onto $T_X\iStkn$ and $T_X\iStkn^{\perp}_{\rho}$ are respectively given by 
\begin{align}
	\PZ(Y) & = Y - XJ\sym(X^TAY)\label{eq:eq37a1}\\
	& = XJ\skew(X^TAY)+(I_n-XJX^TA)Y\label{eq:eq37a2}
\end{align}
and 
\begin{equation}	
	\PZ^{\perp}(Y) = XJ\sym(X^TAY). \label{eq:eq37c}
\end{equation}
\end{proposition}

{\it Proof}
From \eqref{eq:eq31b} and Proposition~\ref{prop:normsp2}, the images of $Y$ should take the forms
$ \PZ(Y) = XW_Y + A^{-1}X_{\perp}K_Y $ and $ \PZ^{\perp}(Y) = XU_Y,$ in which $JW_Y\in\skewset(k),$ $K_Y\in\mathbb{R}^{(n-k)\times k},$ and $ JU_Y\in\symset(k).$ In view of \eqref{eq:orthogonalsum}, it holds that
\begin{equation}
	Y = 
	XW_Y + A^{-1}X_{\perp}K_Y + XU_Y.\label{eq:eq38}
\end{equation}
Left multiplying 
this identity by $X^TA$  yields
$X^TAY = 
JW_Y + JU_Y$.
Therefore, $\sym(X^TAY) = 
JU_Y$ leads to $\PZ^{\perp}(Y) = XU_Y = XJ\sym(X^TAY)$ which proves \eqref{eq:eq37c}.

The identity \eqref{eq:eq37a1} follows from \eqref{eq:eq38},
$$
\PZ(Y) = XW_Y + A^{-1}X_{\perp}K_Y = Y - XU_Y 
=  Y - XJ\sym(X^TAY).
$$
Since $-\sym(X^TAY)=\skew(X^TAY) - X^TAY,$ \eqref{eq:eq37a2} is a consequence of \eqref{eq:eq37a1}.
\qed
\subsection{Riemannian gradient}
We discuss now the Riemannian gradient of the cost function $f$ at $X\in\iStkn$ with respect to the generalized canonical metric $g_{\rho,X_{\perp}}$ constructed in Section~\ref{sec:metric}. 
Recalling that it is the unique vector of $T_X\iStkn$, denoted by $\gradc(X)$,  satisfying
\begin{equation*}
g_{\rho,X_{\perp}}(\gradc(X),Z) = \Df \bar{f}(X)[Z] = 
\langle\nabla\bar{f}(X),Z \rangle
\text{ for all } Z\in T_X\iStkn,
\end{equation*}
where $\bar{f}$ is any smooth extension of $f$ around $X$ in $\Rbnk.$
We have the following formula of $\gradc(X).$ 
\begin{proposition}\label{prop:gradf_gen}
The Riemannian gradient of the cost function $f$ at $X$ in $\iStkn$ is given by
	\begin{equation}
		\gradc(X) = \rho\,XJ\skew(JX^T\nabla\bar{f}(X)) + A^{-1}X_{\perp}\Gamma_3 X_{\perp}^TA^{-1}\nabla\bar{f}(X). \label{eq:eq312a}
	\end{equation}
\end{proposition}

{\it Proof} First, recall that  
the metric $g_{\rho,X_{\perp}}$ in \eqref{eq:canon_metric} can be extended to $\mathbb{R}^{n\times k}$, and denote by $\mathrm{grad}_{gc}\bar{f}(X)$ the Riemannian gradient of $\bar{f}$ on $\mathbb{R}^{n\times k}$ equipped with that extended metric. In view of the definition of the Riemannian gradient,  
there holds
$$ 
\langle\Mx \mathrm{grad}_{gc}\bar{f}(X),Z\rangle 
= g_{\rho,X_{\perp}}\left(\mathrm{grad}_{gc}\bar{f}(X),Z\right) =\Df\bar{f}(X)[Z] = 
\langle\nabla\bar{f}(X),Z \rangle,$$
for any $Z\in \Rbnk.$ Therefore, from the formula  \eqref{eq:explicitMx_2} for the inverse of $\Mx$,
\begin{equation*}
	\mathrm{grad}_{gc}\bar{f}(X) = \Mxa \nabla\bar{f}(X) 
	= \left( \rho XX^T + A^{-1}X_{\perp}\Gamma_3 X_{\perp}^TA^{-1}\right)\nabla\bar{f}(X). 
\end{equation*}
It follows from \cite[Eq.~(3.37)]{AbsiMS08} and \eqref{eq:eq37a2} that
\begin{align*}
	\gradc(X) &= \PZ(\mathrm{grad}_{gc}\bar{f}(X))\\
	&=  XJ\skew\biggr(X^TA\left( \rho XX^T + A^{-1}X_{\perp}\Gamma_3 X_{\perp}^TA^{-1}\right)\nabla\bar{f}(X) \biggr)\\
	&\quad + \;(I_n -XJX^TA)\left( \rho XX^T + A^{-1}X_{\perp}\Gamma_3 X_{\perp}^TA^{-1}\right)\nabla\bar{f}(X)\\
	&= \rho XJ\skew\left(J X^T\nabla\bar{f}(X)\right) +
	\left( \rho XX^T + A^{-1}X_{\perp}\Gamma_3 X_{\perp}^TA^{-1}\right)\nabla\bar{f}(X)\\
	&\quad - \rho XX^T \nabla\bar{f}(X)\\
	&= \rho XJ\skew\left(J X^T\nabla\bar{f}(X)\right) + A^{-1}X_{\perp}\Gamma_3 X_{\perp}^TA^{-1}\nabla\bar{f}(X).
\end{align*}
This proves the identity \eqref{eq:eq312a}.
	\qed
	
	It is stressed that in practice, the appearance of $X_\perp$ in \eqref{eq:eq312a} needs to be removed. This must be made in accordance with the  
	condition imposed in either \eqref{eq:Xperp_1} or \eqref{eq:Xperp_2}.
	Consider first the case $\left(X_\perp^TX_\perp\right)^{-1}=\Gamma_3$. It is readily checked that 
	$$X_\perp\Gamma_3 X_\perp^T = X_\perp\left(X_\perp^TX_\perp\right)^{-1} X_\perp^T = I_n-X(X^TX)^{-1}X^T.$$ 
	Hence, the expression \eqref{eq:eq312a} becomes
	\begin{equation}\label{eq:grad_practice_1}
		\begin{aligned}
			\gradc(X) =& \rho\,XJ\skew(JX^T\nabla\bar{f}(X))\\ 
			&+A^{-1}\left(I_n-X(X^TX)^{-1}X^T\right)A^{-1}\nabla\bar{f}(X).
		\end{aligned}
	\end{equation}
	In the case that $ \left(X_{\perp}^TA^{-1}X_{\perp}\right)^{-1}X_{\perp}^TX_{\perp}\left(X_{\perp}^TA^{-1}X_{\perp}\right)^{-1} = \Gamma_3$, taking the formula \eqref{eq:In_1} resulting from the identity $I_n = EE^{-1}$ into account, it follows that
	\begin{align*}
		A^{-1}X_{\perp}\Gamma_3X_{\perp}^TA^{-1} &= A^{-1}X_{\perp}\left(X_{\perp}^TA^{-1}X_{\perp}\right)^{-1}X_{\perp}^TX_{\perp}\left(X_{\perp}^TA^{-1}X_{\perp}\right)^{-1}X_{\perp}^TA^{-1}\\
		&= \left(I_n-XJX^TA\right)\left(I_n-XJX^TA\right)^T.
	\end{align*}
	Thus, we obtain
	\begin{equation}\label{eq:grad_practice_2}
		\begin{aligned}
			\gradc(X) =& \rho\,XJ\skew(JX^T\nabla\bar{f}(X)) \\
			&+ \left(I_n-XJX^TA\right)\left(I_n-XJX^TA\right)^T\nabla\bar{f}(X).
		\end{aligned}
	\end{equation}
	If matrix $A$ is dense and has no additional special structure, the computational cost for $A^{-1}(\cdot)$ in \eqref{eq:grad_practice_1} is generally comparable to $n^3$. Therefore, in this case, using  \eqref{eq:grad_practice_2} is more advantageous. One can see that the computation of \eqref{eq:grad_practice_2} requires $2n^2k + 10nk^2 - nk -3k^2$ flops. Obviously, when $k$ is 
	of the same scale as $n$, i.e., $k = \mathcal{O}(n)$, this cost is cheaper than that caused by \eqref{eq:gradMx_formula}. This will be seen in the numerical examples presented in Section~\ref{sec:NumerExam}.

	Before closing this section, we discuss the reason for some closed-form expressions of the orthogonal projections in popular instances. Indeed, our general setting here and the tractable metric in \cite[Sect. 3]{TiepSon2025} 
	allow a clear explanation. First, for the orthogonal Stiefel manifold, the canonical metric proposed in \cite[Subsect. 2.3]{EdelAS98} is a special case of our proposed metric \eqref{eq:canon_metric} when $\rho = 2$ and $\Gamma_3=I_{n-k}$. Therefore, there is no surprise when the associated orthogonal projection on the tangent space is in closed-form; namely, \cite[Eq.~(2.4)]{EdelAS98} is a special instance of \eqref{eq:eq37a2}. Alternatively, if the Euclidean metric is used, i.e., $\Mx = I_n$, equation \eqref{eq:Lya_Uxy} immediately admits the explicit solution $U_{X,Y} = \sym(X^TY)$, which together with \cite[Prop. 3.2]{TiepSon2025}, yields \cite[Eq.~(2.4)]{EdelAS98}. Second, for the generalized Stiefel manifold, it was suggested in \cite{ShusA2023,WangDPY2024} to choose $\Mx = A$. This results in the explicit solution $U_{X,Y}= \sym(X^TAY)$ to \eqref{eq:Lya_Uxy}. In view of this, \cite[Prop. 3.2]{TiepSon2025} implies \cite[Eq.~(3.18)]{ShusA2023} and \cite[Eq.~(4.3)]{WangDPY2024}. Note, however, that unlike our derivation, \cite[Eq.~(3.18)]{ShusA2023} was obtained from \cite[Eq.~(2.4)]{EdelAS98} via a change of variable.
	
	In spite of the difference, similar arguments also work for  
	the geometry of the symplectic Stiefel manifold under the canonical-like metric \cite{GSAS21} as a special case of the tractable metric \cite{GSS2024_2}.
	\section{Quasi-geodesics with respect to the generalized canonical metric}\label{sec:quasigeo}
	A geodesic on a manifold is a smooth curve with vanishing acceleration. In the dawn of Riemannian optimization, the solution is updated by searching along the geodesic which goes 
	through the current solution in an appropriate tangent direction \cite{EdelAS98,Gaba82,Smit94}. In the Euclidean space, zero acceleration means 
	$\ddot{X}(t)=0$, 
	or equivalently $X(t) = \dot{X}(0)t+X(0)$, which 
	yields a straight line. That is, geodesics can be viewed as a generalization of lines in the Euclidean space to the Riemannian manifold. In a general setting, the notion of acceleration depends on the Riemannian connection \cite[Sect. 5.4]{AbsiMS08} and therefore the differential equation determining geodesics is usually much more involved. This makes the task of finding an explicit formula for the geodesic, given initial information $X(0), \dot{X}(0)$, very difficult. Fortunately, in optimization, only local information encoded by geodesics and its feasibility are 
	interesting due to their usage in updating the iterate. Thus, one usually retains these features in approximating geodesics which results in the so-called \emph{quasi-geodesics}; see, e.g., \cite{JurdML2020,KrakMLB2017,NishA2005} for more approaches  
	to the quasi-geodesics on the orthogonal Stiefel manifold. 
	
	In this paper, we follow a procedure similar to  that of \cite[Subsect.~4.3]{GSAS21}, which was adopted from \cite[Subsect.~2.2]{EdelAS98}. Namely, we replace the acceleration that should be formulated in accordance with the Riemannian connection associated with the proposed generalized canonical metric by the one 
	corresponding to the Euclidean metric, which turns out to be $\ddot{Y}(t).$ In view of the normal vectors formulated in Proposition~\ref{prop:normsp2}, we define a quasi-geodesic with respect to the  
	generalized canonical metric $g_{\rho,X_{\perp}}$ 
	to be any solution $Y_{\rho}(t)\in\iStkn$ of  
	the initial value problem 
	\begin{equation}\label{eq:quasiwrtgrh}
		\ddot{Y}(t) + Y(t) U(t) =0,\ Y(0) = X\in\iStkn,\ \dot{Y}(0) = Z\in T_X\iStkn,
	\end{equation}
	with some $U(t)$ such that $JU(t)\in\symset(k)$ for all $t\geq 0$, although the obtained result  apparently holds for $t<0$ too. In this section, a closed-form solution to this equation is presented. 
	For the sake of brevity, we refrain from writing the variable $t$ in our presentation except for the arguments in which the dependence on $t$ is highlighted. 
	
	As the first step, we show that for the equation \eqref{eq:quasiwrtgrh} to have a solution, $U(t)$ must have a special form. Consecutively 
	differentiating with respect to $t$ both sides of the identity $Y^TAY = J$ twice
	yield
	\begin{align}\label{eq:deriv_eq_1}
		\dot{Y}^TAY + Y^TA\dot{Y} &=0,\\
		\ddot{Y}^TAY + 2\dot{Y}^TA\dot{Y} + Y^TA\ddot{Y} & =0,\label{eq:deriv_eq_2}
	\end{align}
	respectively. Substituting $\ddot{Y} = -YU$ from \eqref{eq:quasiwrtgrh} into \eqref{eq:deriv_eq_2} gives
	\begin{equation*}
		(-YU)^TAY + 2\dot{Y}^TA\dot{Y} + Y^TA(-YU)=0 
	\end{equation*}
	or $\dot{Y}^TA\dot{Y}=JU$ due to the symmetry of $JU$. Hence, $U = J\dot{Y}^TA\dot{Y}$ and therefore, 
	the equation \eqref{eq:quasiwrtgrh} becomes
	\begin{equation}\label{eq:quasigrho2}
		\ddot{Y} + YJ(\dot{Y}^TA\dot{Y})=0.
	\end{equation}
	
	Next, let us set 
	$S \coloneqq Y^TA\dot{Y}$ and $V\coloneqq\dot{Y}^TA\dot{Y}$. Thanks to the relation \eqref{eq:deriv_eq_1}, 
	$S\in\skewset(k)$ for $t\geq 0.$  Moreover, using \eqref{eq:quasigrho2}, we get
	$$\dot{S} = V + Y^TA\ddot{Y} = V - Y^TAYJV = V-V=0, \mbox{ for } t\geq 0.$$
	Hence, $S(t) = S(0) \eqqcolon S_0\in
	\skewset(k)$ for $t\geq 0$. We proceed to determine $V$ which is symmetric. 
	Taking the derivative of $V$ with \eqref{eq:quasigrho2} in mind, we have that
	\begin{align*}
		\dot{V} &= \ddot{Y}^TA\dot{Y} + \dot{Y}^TA\ddot{Y} = (-YJV)^TA\dot{Y} + \dot{Y}^TA(-YJV)\\
		&= -VJS + SJV = S_0JV-VJS_0. 
	\end{align*}
	We will show that the initial value problem 
	\begin{equation}\label{eq:quasigrho3}
		\dot{V} = S_0JV - VJS_0,\ V(0) = V_0,
	\end{equation}
	admits a solution $V(t) = \e^{-tS_0^TJ}V_0\e^{-tJS_0}.$ 
	Indeed, differentiating such $V(t)$ yields
	\begin{align*}
		\dot{V} 
		& = -S_0^TJ\e^{-tS_0^TJ} V_0 \e^{-tJS_0} -\e^{-tS_0^TJ} V_0 \e^{-tJS_0}JS_0\\
		& = -S_0^TJV -VJS_0 = S_0JV -VJS_0.
	\end{align*}
	
	Next, we present a technical lemma, which is similar to \cite[Prop. 4.6]{GSAS21}. 
	\begin{lemma}\label{lem:invariance_Exp}
		If $\Theta\in\skewset(k)$ and $\Gamma\in\symset(k),$ then for any $t\in \Rb,$ it holds that
		\begin{equation}\label{eq:invariance_Exp}
			\left(\exp\left(t\begin{bmatrix}
				J\Theta& J\Gamma\\
				I_k&J\Theta
			\end{bmatrix}\right)\right)^T\begin{bmatrix}
				0& J\\
				-J&0
			\end{bmatrix}\exp\left(t\begin{bmatrix}
				J\Theta& J\Gamma\\
				I_k&J\Theta
			\end{bmatrix}\right)=\begin{bmatrix}
				0& J\\
				-J&0
			\end{bmatrix}.\end{equation}
	\end{lemma}
	
{\it Proof} For the sake of convenience, let us set $\Lambda\coloneq \left[\begin{smallmatrix}
			J\Theta& J\Gamma\\
			I_k&J\Theta
		\end{smallmatrix}\right]$ and the function on the left hand side of \eqref{eq:invariance_Exp} as $u(t)\coloneq \left(\exp\left(t\Lambda\right)\right)^T\left[\begin{smallmatrix}
			0& J\\
			-J&0
		\end{smallmatrix}\right]\exp\left(t\Lambda\right)$. Then, it follows that 
		$$\dot{u}(t) = \left(\exp\left(t\Lambda\right)\right)^T\biggr(\Lambda^T\begin{bmatrix}
			0& J\\
			-J&0
		\end{bmatrix}+\begin{bmatrix}
			0& J\\
			-J&0
		\end{bmatrix}\Lambda\biggr)\exp\left(t\Lambda\right).$$ 
		By assumption, one can easily check that 
		$$\Lambda^T\begin{bmatrix}
			0& J\\
			-J&0
		\end{bmatrix}+\begin{bmatrix}
			0& J\\
			-J&0
		\end{bmatrix}\Lambda = 0.$$
		Therefore, $\dot{u}(t) = 0$ which leads to the fact that $u(t) = u(0) =\begin{bmatrix}
			0& J\\
			-J&0
		\end{bmatrix},$ for any $t\in\Rb.$ 
	\qed
	
	Now, we are ready to give the explicit formula for a solution  
	to the problem  \eqref{eq:quasiwrtgrh}. 
	\begin{proposition}\label{prop:expl_sol}
		Given $X\in\iStkn,Z\in T_X\iStkn.$ The initial value problem
		\begin{equation}\label{eq:quasiwrtgrh_2}
			\ddot{Y}(t) + Y(t)J\dot{Y}(t)^TA\dot{Y}(t) =0,\ Y(0) = X,\ \dot{Y}(0) = Z ,
		\end{equation}
		has a unique global solution in $\iStkn$ given by
		\begin{equation}\label{eq:Yrho}
			Y_{\rho}(t) = Y_{\rho}^{\mathrm{qgeo}}(t,Z)\coloneqq [X \; \; Z]\exp\left( t\begin{bmatrix}
				JS_0& -JV_0\\
				I_k&JS_0
			\end{bmatrix}\right)\begin{bmatrix}
				I_k\\
				0
			\end{bmatrix}\e^{-tJS_0},
		\end{equation}
		where $S_0 
		 = Y_{\rho}^T(0)A\dot{Y}_{\rho}(0) =  X^TAZ$ and $ V_0 \coloneq \dot{Y}_{\rho}^T(0)A\dot{Y}_{\rho}(0) =  Z^TAZ.$
	\end{proposition}
	
	{\it Proof}
		First, let us assume that the problem \eqref{eq:quasiwrtgrh_2} has a global solution in $\iStkn,$ denoted also by $Y_{\rho}(t)$. It follows that 
		\begin{align*}
			\ddt\left(Y_{\rho}(t) \e^{tJS_0}\right) & =\dot{Y_{\rho}}(t)\e^{tJS_0}+Y_{\rho}(t)\e^{tJS_0}JS_0.\\
			\ddt\left(\dot{Y}_{\rho}(t) \e^{tJS_0}\right) & =\ddot{Y_{\rho}}(t)\e^{tJS_0}+\dot{Y}_{\rho}(t)\e^{tJS_0}JS_0\\
			&= \left( -Y_{\rho}(t)JV\right) \e^{tJS_0} + \dot{Y}_{\rho}(t)\e^{tJS_0}JS_0\\
			& =-Y_{\rho}(t)J\e^{-tS_0^TJ}V_0\e^{-tJS_0} \e^{tJS_0}+ \dot{Y}_{\rho}(t)\e^{tJS_0}JS_0\\
			&= 	-Y_{\rho}(t)\e^{tJS_0}JV_0+\dot{Y}_{\rho}(t)\e^{tJS_0}JS_0,
		\end{align*}
		in which we have invoked \eqref{eq:quasigrho3} to replace 
		$V \coloneq \dot{Y_{\rho}}^TA\dot{Y_{\rho}}$ by its explicit form 
		$V =\e^{-tS_0^TJ}V_0\e^{-tJS_0}$ 
		when $Y_{\rho}(t)\in\iStkn,$ for all $t\geq 0.$ 
		Therefore,
		\begin{equation}\label{eq:couple_eq}
			\ddt\left[Y_{\rho}(t) \e^{tJS_0} \;\; \dot{Y}_{\rho}(t) \e^{tJS_0}\right] =  \left[Y_{\rho}(t) \e^{tJS_0} \;\; \dot{Y}_{\rho}(t) \e^{tJS_0}\right]\begin{bmatrix}
				JS_0&-JV_0\\
				I_k&JS_0
			\end{bmatrix}.
		\end{equation}
		This basic equation admits an explicit solution from which we obtain that
		\begin{align*}
			Y_{\rho}(t) &= \left[Y_{\rho}(0) \;\; \dot{Y}_{\rho}(0)\right]\exp\left( t\begin{bmatrix}
				JS_0& -JV_0\\
				I_k&JS_0
			\end{bmatrix}\right)\begin{bmatrix}
				I_k\\
				0
			\end{bmatrix}\e^{-tJS_0}\\
			&=[X\;\; Z]\exp\left( t\begin{bmatrix}
				JS_0& -JV_0\\
				I_k&JS_0
			\end{bmatrix}\right)\begin{bmatrix}
				I_k\\
				0
			\end{bmatrix}\e^{-tJS_0},
		\end{align*}
		which is the expected formula.
		
		Next, we will show that this is indeed a global solution in $\iSt(k,n)$ to the problem \eqref{eq:quasiwrtgrh_2}. For future reference, we present the feasibility part in Proposition~\ref{prop:global_feasible} below. In the interest of brevity, 
		we set $\Psi \coloneq \left[\begin{smallmatrix}
			JS_0& -JV_0\\
			I_k&JS_0
		\end{smallmatrix}\right]$. In view of  \eqref{eq:Yrho}, it holds that
		\begin{equation*}
			\dot{Y}_{\rho}\e^{tJS_0}+Y_{\rho}\e^{tJS_0}JS_0 = \ddt \left(Y_{\rho}\e^{tJS_0}\right) = [X \; \; Z]\exp\left( t\Psi\right)\Psi\begin{bmatrix}
				I_k\\
				0
			\end{bmatrix}.
		\end{equation*}
			This leads to
			\begin{equation}\label{Yrho_dotExp}
				\begin{aligned}
					\dot{Y}_{\rho}\e^{tJS_0} =& [X \; \; Z]\exp\left( t\Psi\right)\Psi\begin{bmatrix}
						I_k\\
						0
					\end{bmatrix}- [X \; \; Z]\exp\left( t\Psi\right)\begin{bmatrix}
						I_k\\
						0
					\end{bmatrix}JS_0 \\
					=& [X \; \; Z]\exp\left( t\Psi\right)\biggr(\begin{bmatrix}
						JS_0\\
						I_k
					\end{bmatrix} -\begin{bmatrix}
						JS_0\\
						0
					\end{bmatrix}\biggr) = [X \; \; Z]\exp\left( t\Psi\right)\begin{bmatrix}
						0\\
						I_k
					\end{bmatrix}.
				\end{aligned}
			\end{equation}
			Repeating this procedure with $\dot{Y}_{\rho}\e^{tJS_0}$, we obtain
				\[ \ddot{Y}_{\rho}\e^{tJS_0}+\dot{Y}_{\rho}\e^{tJS_0}JS_0 = [X \; \; Z]\exp\left( t\Psi\right)\Psi\begin{bmatrix}
					0\\
					I_k
				\end{bmatrix}.\]
				Using \eqref{Yrho_dotExp}, this identity gives 
				\begin{equation}\label{LHS_quasi_geo_ode}
					\ddot{Y}_{\rho}\e^{tJS_0} = [X \; \; Z]\exp\left( t\Psi\right)\begin{bmatrix}
						-JV_0\\
						0
					\end{bmatrix}.
				\end{equation}
				On the other hand, using the expression for $\dot{Y}_{\rho}$ in \eqref{Yrho_dotExp}, we receive the identities 
				\begin{align}\notag
					\dot{Y}_{\rho}^TA\dot{Y}_{\rho} =& \left(\e^{-tJS_0}\right)^T [0 \; \; I_k]\left(\exp(t\Psi)\right)^T\begin{bmatrix}
						X^T\\
						Z^T
					\end{bmatrix}A[X \; \; Z]\exp\left( t\Psi\right)\begin{bmatrix}
						0\\
						I_k
					\end{bmatrix}\e^{-tJS_0}\\ \notag
					=& \,\e^{-tS_0^TJ} [0 \; \; I_k]\left(\exp(t\Psi)\right)^T\begin{bmatrix}
						J&S_0 \\
						-S_0&V_0
					\end{bmatrix}\exp\left( t\Psi\right)\begin{bmatrix}
						0\\
						I_k
					\end{bmatrix}\e^{-tJS_0}\\ \notag
					=&\, \e^{-tS_0^TJ}[0 \; \; I_k]\left(\exp(t\Psi)\right)^T\begin{bmatrix}
						0&J\\
						-J&0
					\end{bmatrix}\Psi\exp\left( t\Psi\right)\begin{bmatrix}
						0\\
						I_k
					\end{bmatrix}\e^{-tJS_0}\\ \notag
					=&  \,\e^{-tS_0^TJ} [0 \; \; I_k]\left(\exp(t\Psi)\right)^T\begin{bmatrix}
						0&J\\
						-J&0
					\end{bmatrix}\exp\left( t\Psi\right)\Psi\begin{bmatrix}
						0\\
						I_k
					\end{bmatrix}\e^{-tJS_0}\\ \notag
					=& \,\e^{-tS_0^TJ}[0 \; \; I_k]\begin{bmatrix}
						0&J\\
						-J&0
					\end{bmatrix} \Psi\begin{bmatrix}
						0\\
						I_k
					\end{bmatrix}\e^{-tJS_0}\\
					=&\,\e^{-tS_0^TJ}V_0 \e^{-tJS_0}
					\label{Yrho_dotT_A_Yrho_dot},
				\end{align}
				in which Lemma \ref{lem:invariance_Exp} have been invoked in the fifth identity.
						
						At the end, left multiplying 
						\eqref{Yrho_dotT_A_Yrho_dot} by $Y_{\rho}J$ yields
							\begin{align}\notag
								Y_{\rho}J(\dot{Y}_{\rho}^TA\dot{Y}_{\rho})=&\,Y_{\rho}J\e^{-tS_0^TJ}V_0 \e^{-tJS_0}\\ \notag
								=&\,[X \; \; Z]\exp\left( t\Psi\right)\begin{bmatrix}
									I_k\\
									0
								\end{bmatrix}\e^{-tJS_0}J^2\e^{tJS_0}JV_0\e^{-tJS_0}\\ \notag
								=&\, [X \; \; Z]\exp\left( t\Psi\right)\begin{bmatrix}
									I_k\\
									0
								\end{bmatrix}JV_0 \e^{-tJS_0}\\
								= &\, [X \; \; Z]\exp\left( t\Psi\right)\begin{bmatrix}
									JV_0\\
									0
								\end{bmatrix}\e^{-tJS_0}, \label{RHS_quasi_geo_ode}
							\end{align}
						in which the second identity is due to 
						the fact that
						\begin{equation}\label{eq:extra1_proof}
							 (\e^{-tJS_0})^T=\e^{-tS_0^TJ} = J^2\e^{-tS_0^TJ} = J\e^{tJS_0}J.
						\end{equation}
						The conclusion follows from the two identities \eqref{LHS_quasi_geo_ode} and \eqref{RHS_quasi_geo_ode}. 
						
						The uniqueness of the solution is derived from that of the problems \eqref{eq:quasigrho3} and \eqref{eq:couple_eq} which is a well-known 
						fact in the theory of linear differential equations in Banach spaces; see, e.g., \cite[Sect. 1.1]{Deimling1977}.
					
					\begin{proposition}\label{prop:global_feasible}
						The curve $Y_\rho(t)$ given in \eqref{eq:Yrho} lies in the indefinite Stiefel manifold $\iSt_{A,J}(k,n)$, i.e., for any $t\geq 0,$ we have $Y_{\rho}(t)^TAY_{\rho}(t)=J.$
					\end{proposition}
					
					{\it Proof}
						From \eqref{eq:Yrho}, it follows that 
						\begin{equation}\label{YrhoT_A_Yrho}
							Y_{\rho}^TAY_{\rho} =  \left(\e^{-tJS_0}\right)^T [I_k \; \; 0]\left(\exp(t\Psi)\right)^T\begin{bmatrix}
								X^T\\
								Z^T
							\end{bmatrix}A[X \; \; Z]\exp\left( t\Psi\right)\begin{bmatrix}
								I_k\\
								0
							\end{bmatrix}\e^{-tJS_0}.
						\end{equation}
						According to the proof 
						for \eqref{Yrho_dotT_A_Yrho_dot},  
						we have that
						\begin{equation*}
							\left(\exp(t\Psi)\right)^T\begin{bmatrix}
								X^T\\
								Z^T
							\end{bmatrix}A[X \; \; Z]\exp\left( t\Psi\right)
							=  \begin{bmatrix}
								0&J \\
								-J&0
							\end{bmatrix}\Psi = \begin{bmatrix}
								J&S_0\\
								-S_0&V_0
							\end{bmatrix}.
						\end{equation*}
						This fact leads \eqref{YrhoT_A_Yrho} to 
						\begin{align*}
							Y_{\rho}^TAY_{\rho} =&  \left(\e^{-tJS_0}\right)^T [I_k \; \; 0]\begin{bmatrix}
								J&S_0\\
								-S_0&V_0
							\end{bmatrix}
							\begin{bmatrix}
								I_k\\
								0
							\end{bmatrix}\e^{-tJS_0}\\
							=&J\e^{tJS_0}J^2\e^{-tJS_0}=J,
						\end{align*}
						in which the second identity follows from the fact that $\left(\e^{-tJS_0}\right)^T = J\e^{tJS_0}J$ derived in \eqref{eq:extra1_proof}. 
					\qed

					In the case of generalized Stiefel manifold, i.e.,  $A$ is spd and $J=I_k$, the obtained quasi-geodesic reduces to a simple form as follows.
					\begin{corollary}
						Given $X\in\St_A(k,n)$, $Z\in T_X\St_A(k,n)$, and set $S_0 =  X^TAZ$, $ V_0 =  Z^TAZ$. The curve 
						\begin{equation}\label{eq:qgeo_genStiefel}
							Y(t) = [X \; \; Z]\exp\left( t\begin{bmatrix}
								S_0& -V_0\\
								I_k&S_0
							\end{bmatrix}\right)\begin{bmatrix}
								I_k\\
								0
							\end{bmatrix}\e^{-tS_0}
						\end{equation}
						is a quasi-geodesic associated with $X$ and $Z$ in the sense that $Y(t) \in \St_A(k,n)$ for all $t\geq 0$ and $Y(0) = X, \dot{Y}(0) = Z$.
					\end{corollary}
					\begin{remark}
						Restricting to the orthogonal Stiefel manifold,  the quasi-geodesic \eqref{eq:qgeo_genStiefel} 
						reduces to the  geodesic obtained in \cite[Subsect. 2.2]{EdelAS98}. It is worth 
						noting that the one in \cite[Subsect. 2.2]{EdelAS98} is indeed associated with the  Euclidean metric. Nevertheless, by  direct verification, 
						one can show that 
						the normal space, orthogonal projections onto the tangent and normal spaces under the Euclidean and under the canonical metrics are identical. Therefore, for the orthogonal Stiefel manifold, the embedded geodesic associated with the Euclidean metric is also the quasi-geodesic associated with the canonical metric. Hence, \eqref{eq:Yrho} is a generalization of this quasi-geodesic to the indefinite Stiefel manifold.
					\end{remark}
					
					
					\section{Quasi-geodesic based retraction on $\iStkn$}\label{sec:rectraction}
					Let us first formally define a retraction on $\iStkn$; see \cite{AbsiMS08,AdleDMMS02} for more  details about retractions on a general Riemannian manifold. It is a smooth mapping $\RX$ from the tangent bundle $T\iStkn\coloneqq \bigcup_{X\in\iStkn}T_X\iStkn$ to $\iStkn$ satisfying the following conditions for all $X\in\iStkn:$
					\begin{align*}
						\tag{R.1} \; \;\; &\RZ(0) = X, \text{ where } 0 \text{ is the origin of } T_X\iStkn,\label{eq:conditionC1}\\
						\tag{R.2} \; \; \;& \ddt\RZ\left(tZ\right)\rvert_{t=0} = Z, \text{ for all } Z\in T_X\iStkn,\label{eq:conditionC2}
					\end{align*}
					where $\RZ$ is the
					restriction of $\RX$ to $T_X\iStkn.$ The crucial role of a retraction in optimization is to maintain the feasibility of the updated iterate from a feasible search direction which is a tangent vector at the current iterate. In \cite{TiepSon2025}, a retraction based on the Cayley transformation was proposed. In this work, based on the formulation of a quasi-geodesic \eqref{eq:Yrho}, 
					we propose another retraction, which  will be termed \emph{quasi-geodesic retraction}, as follows:
					\begin{align}%
						\notag
						\RX^{\mathrm{qgeo}}: T\iStkn & \longrightarrow \iStkn & &\\ \label{eq:retractiongrho}
						Z\in T_X\iStkn &\longmapsto \RZ^{\mathrm{qgeo}}(Z) \coloneqq  Y_{\rho}^{\mathrm{qgeo}}(1,Z) &&\\ 
						&= [X\;\; Z]\exp\left( \begin{bmatrix}
							JX^TAZ& -JZ^TAZ\\
							I_k&JX^TAZ
						\end{bmatrix}\right)\begin{bmatrix}
							I_k\\
							0
						\end{bmatrix}\e^{-JX^TAZ},&&
						\notag
					\end{align}
					where $\RZ^{\mathrm{qgeo}}$ is the restriction of $\RX^{\mathrm{qgeo}}$ to $T_X\iStkn.$ The well-definedness of $\RX^{\mathrm{qgeo}}$ is guaranteed in the following result whose proof is similar to \cite[Lem. 5.1]{GSAS21}.
					\begin{proposition}
						The map $\RX^{\mathrm{qgeo}}$ given in \eqref{eq:retractiongrho} 
						is a globally defined retraction on $\iStkn.$
					\end{proposition}
					
					{\it Proof} We only need to verify the conditions \eqref{eq:conditionC1} and \eqref{eq:conditionC2} since the expression in \eqref{eq:retractiongrho} 
						is well-defined for any $X\in\iStkn$ and $Z\in T_X\iStkn$ thanks to Proposition~\ref{prop:global_feasible}
						
						It holds that $\RZ^{\mathrm{qgeo}}(0) = [X\;\; 0]\exp\left(\left[\begin{smallmatrix}
							0& 0\\
							I_k&0
						\end{smallmatrix}\right]\right)\left[\begin{smallmatrix}
							I_k\\
							0
						\end{smallmatrix}\right]\e^{0}=[X\;\; 0]\left[\begin{smallmatrix}
							I_k\\
							I_k
						\end{smallmatrix}\right] =X$. 
						Moreover, for any number $t\geq 0,$ 
						\begin{center}
							$\RZ^{\mathrm{qgeo}}(tZ) =
							[X\;\; tZ]\exp\left( \begin{bmatrix}
								tJX^TAZ& -t^2JZ^TAZ\\
								I_k&tJX^TAZ
							\end{bmatrix}\right)\begin{bmatrix}
								I_k\\
								0
							\end{bmatrix}\e^{-tJX^TAZ}.$
						\end{center} 
						So,
						\begin{align*}
							\ddt \RZ^{\mathrm{qgeo}}(tZ) \restrict{t=0}&=
							[0\;\; Z]\exp\left( \begin{bmatrix}
								0& 0\\
								I_k&0
							\end{bmatrix}\right)\begin{bmatrix}
								I_k\\
								0
							\end{bmatrix}\e^{0}\\
							&+[X\;\; 0]\exp\left( \begin{bmatrix}
								0& 0\\
								I_k&0
							\end{bmatrix}\right) \begin{bmatrix}
								JX^TAZ& 0\\
								0&JX^TAZ
							\end{bmatrix}\begin{bmatrix}
								I_k\\
								0
							\end{bmatrix}\e^{0}\\
							&-[X\;\; 0]\exp\left( \begin{bmatrix}
								0& 0\\
								I_k&0
							\end{bmatrix}\right)\begin{bmatrix}
								I_k\\
								0
							\end{bmatrix}\e^{0}JX^TAZ\\
							& = Z,
						\end{align*} 
						where we have invoked the identity $\frac{\d t}{t}\e^{M(t)}|_{t=0} = \e^{M(0)}\dot{M}(0)$ for smooth one-parameter matrix $M(t)$ such that $M(0)\dot{M}(0) = \dot{M}(0)M(0)$; see \cite[Eq.~(2.1)]{Wilcox1967}.
					\qed
					
					\section{Numerical examples}\label{sec:NumerExam} 
					In this section, several numerical examples are presented to verify our theoretical findings. As constructed in Section~\ref{sec:metric}, two formulae for computing the Riemannian gradient using the proposed metric are given in \eqref{eq:grad_practice_1} and \eqref{eq:grad_practice_2}; it was also stated that the latter  outperforms the former when $A$ is dense. First, we validate this statement via an extra example. After that, the derived components are assembled to enable a Riemannian gradient descent scheme and then numerically verified via several model tests. The tests are done on a laptop with Intel(R) Core(TM) i7-4500U (at 1.8GHz, 5MB Cache, 8GB RAM) running MATLAB R2023a under Windows 10 Home.  (The main code is made available at \href{https://sites.google.com/view/ntson/code}{https://sites.google.com/view/ntson/code}.)
						
						\paragraph{A numerical comparison of \eqref{eq:grad_practice_1} and \eqref{eq:grad_practice_2}} We construct a manifold $\iSt_{A,J}$ as follows: the matrix $A = P\diag(1,2,\ldots,1500,-500,-499,\ldots,-1)P^T$, where $P$ is an orthogonalized random matrix 
							and   $J=\diag(I_{14},-I_6)$. We take $f(X)=\frac{1}{2}\tr\left(X^TMX\right)$ with $M$ is a fixed spd matrix of apt size. We compute the Riemannian gradient using two formulae \eqref{eq:grad_practice_1} and \eqref{eq:grad_practice_2} at 
							$X=\left[p_1/\sqrt{1}, \ldots p_{14}/\sqrt{14},p_{1501}/\sqrt{500},\ldots,p_{1506}/\sqrt{495}\right]$ in which $p_j$ is the $j$-th column of $P$ generated previously. We run this setting 10 times and take the average time elapsed. The result shows that, given the Euclidean gradient $\nabla f(X) = MX$, while computing the Riemannian gradient following \eqref{eq:grad_practice_1} takes about 0.45(s),  the one invoking \eqref{eq:grad_practice_2} needs only about 0.016(s), which is roughly 30 times faster. Therefore, only  \eqref{eq:grad_practice_2} is further used in our tests.
					
					\paragraph{A Riemannian gradient descent method} One option to minimize a continuously differentiable cost function $f$ on the indefinite Stiefel manifold $\iSt_{A,J}$,
					\begin{equation}\label{eq:opt_prob}
						\min_{X\in \Rb^{n\times k}} f(X)\ \mbox{ s.t. }\ X^{T}AX = J,
					\end{equation}
					given an initial point $X_0\in \iSt_{A,J}(k,n)$, 
					is to consecutively updates the approximate solution by
					$$X_{j+1} = \RX_{X_j}(-\tau_j \grad f(X_j)),$$ 
					where $\RX_{X_j}(\cdot)$ is a 
					retraction at $X_j$, $\grad f(X_j)$ denotes the Riemannian gradient of the cost function $f$ at $X_j$ with respect to a chosen 
					Riemannian metric $\mathrm{g}(\cdot,\cdot)$, and $\tau_j >0$ is a step size 
					to be determined. 
					To compute the step size, we combine a nonmonotone line search technique with the alternating Barzilai-Borwein strategy 
					on an embedded submanifold; see \cite{BarB88,HuLWY2020,IannP18,ZhanH2004}.  
					For the sake of convenience, we recall the details in  Algorithm~\ref{alg:non-monotone gradient}. 
					Note from \eqref{eq:nonmonotone_cond} that if $\alpha$ is chosen to be zero, this scheme reduces to the popular monotonically decreasing line search with Amijo's rule. Algorithm~\ref{alg:non-monotone gradient} was also employed in \cite{TiepSon2025} for numerical verification of optimization using the geometric quantities under a general tractable metric and the Cayley retraction. Its global convergence is guaranteed by \cite[Thm. 5.7]{GSAS21}. Namely, regardless of the initial point $X_0$, each accumulation point $X^*$ of the iterate generated by Algorithm~\ref{alg:non-monotone gradient} satisfies the first-order necessary optimality condition $\grad f(X^*)=0$. 
					
					\begin{algorithm}[htbp]
						\caption{Riemannian gradient algorithm for the optimization problem~\eqref{eq:opt_prob}}
						\label{alg:non-monotone gradient}
						\begin{algorithmic}[1]
							\REQUIRE The cost function $f$, 
							initial guess $X_0\in\iSt_{A,J}$, $\gamma_0>0$,
							\mbox{$0<\gamma_{\min}<\gamma_{\max}$}, 
							$\beta,$ $\delta\in(0,1)$, $\alpha \in [0,1]$, 
							$q_0=1$, $c_0 = f(X_0)$. 
							\ENSURE Sequence of iterates  $\{X_j\}$. 
							\FOR{$j=0,1,2,\dots$}
							\STATE Compute $Z_j = -\grad f(X_j)$. 
							\IF{$j>0$} 
							\STATE{ 
								$
								\gamma_j=\left\{\begin{array}{ll}
									\dfrac{\langle W_{j-1},W_{j-1}\rangle}{\abs{\tr(W_{j-1}^T Y_{j-1}^{})}} 
									&\text{for odd } j, \\[5mm]
									\dfrac{\abs{\tr(W_{j-1}^T Y_{j-1}^{})}}{\langle Y_{j-1},Y_{j-1}\rangle}
									&\text{for even } j,
								\end{array}\right.
								$\\
								where $W_{j-1} = X_j - X_{j-1}$ and $Y_{j-1} =Z_j-Z_{j-1}$.
							}
							\ENDIF
							\STATE Calculate the trial step size $\gamma_j=\max\bigl(\gamma_{\min},\min(\gamma_j,\gamma_{\max})\bigr)$.			
							\STATE Find the smallest integer $\ell$  such that the 
							nonmonotone condition 
							\begin{equation}\label{eq:nonmonotone_cond}
								f\big(\RX(
								\tau_j Z_j)\big) \le c_j  + \beta\, \tau_j\, \mathrm{g}\big(\grad f(X_j), Z_j\big)
							\end{equation} 
							holds, where $\tau_j=\gamma_j\, \delta^{\ell}$. 
							\STATE Set $X_{j+1} = \RX(\tau_j Z_j)$.
							\STATE Update $q_{j+1} = \alpha q_{j} + 1$ and
							$\displaystyle{c_{j+1} = \frac{\alpha q_{j}}{q_{j+1}} c_{j}  + \frac{1}{q_{j+1}} f(X_{j+1})}$.
							\ENDFOR
						\end{algorithmic}
					\end{algorithm}
					
					\paragraph{Algorithmic parameter setting} 
					For the backtracking search, we set $\beta = 1\e{-4}$, $\delta = 5\e{-1}$, $\gamma_0 = 1\e{-3}$, 
					$\gamma_{\min} = 1\e{-15}$, $\gamma_{\max} =1\e{+5}$.  Moreover, $\alpha = 0.85$ is used in the 
					nonmonotone condition. 
					A running is considered 
					convergent when the norm of the Riemannian gradient of the cost function at the current iterate is relatively small compared to that at the initial guess $X_0$. Namely, $\|\grad f(X_j)\|_{X_j}\leq\texttt{rstop}\|\grad f(X_0)\|_{X_0}$, where $\|\cdot\|_{X_j}$ is the norm associated with the chosen metric at $X_j$. In all tests, \texttt{rstop} is set to be $1\e-5$. We also stop the loop when $j$ exceeds a given maximal number of iterations \texttt{iter} which is set to $1000$.		
					
					\paragraph{Model problems} 
					We examine two 
					problems similar to that in \cite{TiepSon2025}. For the first one, we consider minimizing a trace cost function
					\begin{equation}\label{eq:trace_fun}
						f(X):=\tr(X^TMX),\  X\in \iSt_{A,J}(k,n),
					\end{equation}
					where $n=1000, k = 500$, $M = VV^T$ with $V$ is the $n\times (n-5)$ matrix randomly generated with default setting and then  orthogonalized using \texttt{orth} in MATLAB; matrix $M$ is finally normalized by dividing $M$ by its Frobenius norm. For the manifold, let $p=750, m = 250, k_p = 400, k_m = 100$, we set $A = \diag(1,\ldots,p,-1,\ldots,-m)$ and $J=\diag(I_{k_p},-I_{k_m})$. 
					The initial guess $X_0$ is, for the sake of simplicity, set as the $n\times k$ zero matrix except for the submatrices 
					$X_0(1,\ldots,k_p; 1,\ldots,k_p)= \diag(1/\sqrt{1},\ldots,1/\sqrt{k_p})$ and $X_0(p+1,\ldots,p+k_m; k_p+1,\ldots,k_p+k_m)= \diag( 1/\sqrt{1},\ldots,1/\sqrt{k_m})$. 
					
					The second model is the well-known Procrustes problem
					\begin{equation}\label{eq:procrustes_fun}
						f(X)= \|GX-B\|^2_F
					\end{equation}
					on the $J$-orthogonal group \cite{High03} with $n=500$, $J=\diag(I_p,-I_{m}),$ where $p=375,$ $m = 125$. Numerical data are generated as follows: $G$ is a $(n-20)\times n$ random matrix with default setting and then normalized as above, and $B$ is set as $G\diag(U,V)$ where $U\in \Rb^{p\times p}$ and $V\in\Rb^{m\times m}$ are 
					two orthogonalized random matrices generated with the setting \texttt{rng\,(1,'twister')}. The $n\times n$ identity matrix is used as the initial guess. 
					
					\paragraph{Tests} It is worth
					noting that, while the Riemannian gradient is defined upon the metric, the retraction is metric-independent. 
					Therefore, in Algorithm~\ref{alg:non-monotone gradient}, one can combine any metric with any available retraction. Namely, there are the Cayley retraction and the 
					one based on the quasi-geodesic developed in this paper. Regarding metrics, we have the general tractable and the generalized canonical metrics. In \cite{TiepSon2025}, the Euclidean and weighted Euclidean metrics have been used. The examples there showed that using the latter is much more advantageous. However, 
					the Euclidean Hessian of the cost functions considered here are both degenerate;  
					the choice of weighted Euclidean metric as in \cite{TiepSon2025} is inappropriate. 
					Thus, we will use in our tests the Euclidean and the generalized canonical metrics. Consequently, there are totally four combinations whose corresponding numerical results will be marked by ``Eucl-Cayley'', ``Eucl-qgeo'', ``gcan-Cayley'', and ``gcan-qgeo''.
					
\begin{figure}[ht]
				\centering
				\includegraphics[width=\textwidth]{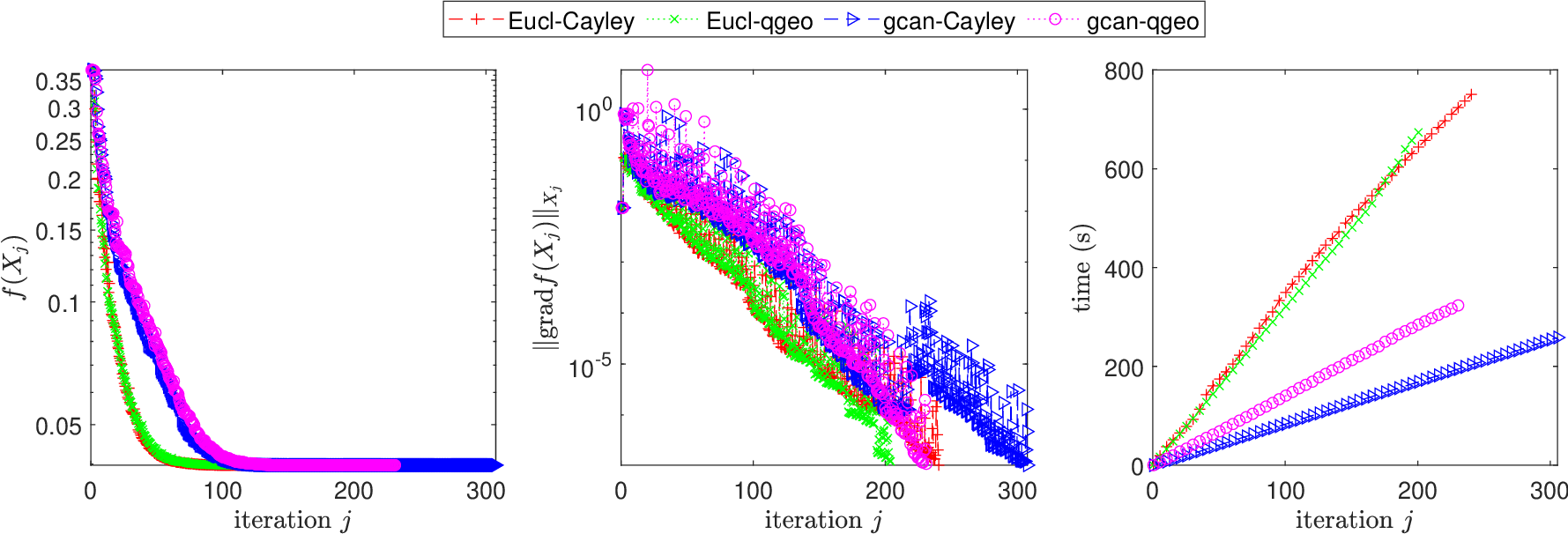}
				\caption{Convergence history of the Riemannian gradient descent method with different combinations of metrics and retractions applied to the trace cost function \eqref{eq:trace_fun}.}\label{fig:tracemin}
			\end{figure}

\begin{figure}[ht]
				\centering
				\includegraphics[width=\textwidth]{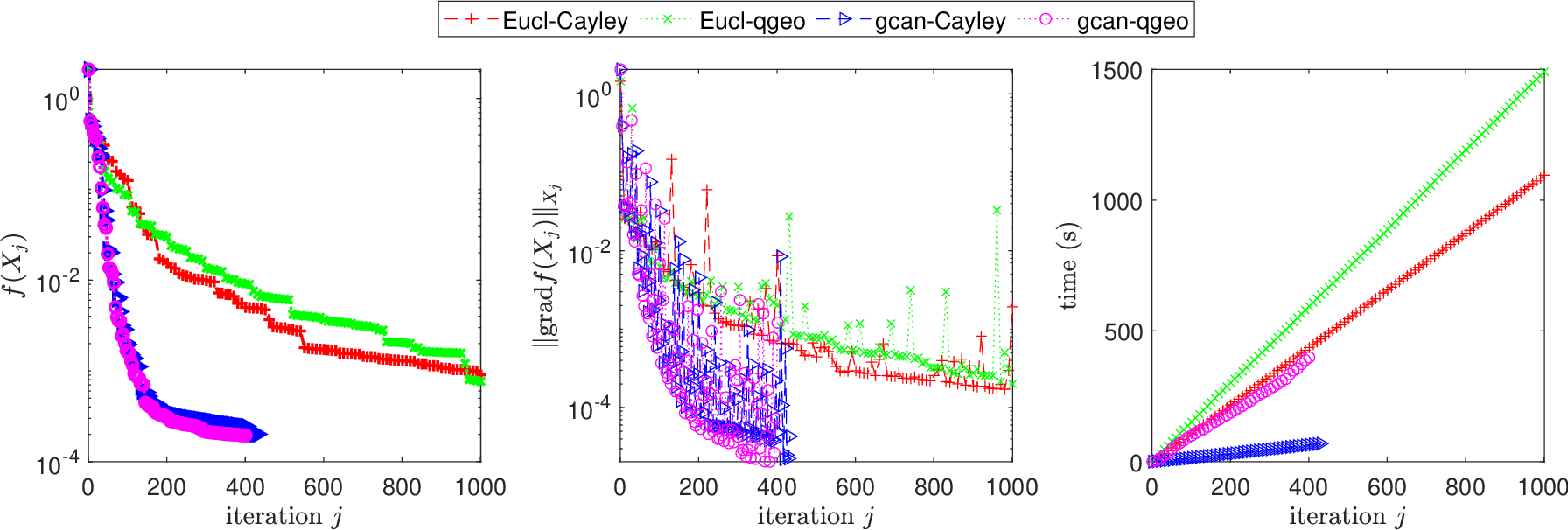}
				{\caption{Convergence history of the Riemannian gradient descent method with different combinations of metrics and retractions applied to the Procrustes cost function \eqref{eq:procrustes_fun}.}\label{fig:procrustes}}
			\end{figure}

					We report the numerical results in Table~\ref{tab:tracemin} and Table~\ref{tab:procrustes_500} in which obj., grad., feas., $\#$iter., $\#$eval., and CPU represent the value of the cost function, the norm in the chosen metric of the Riemannian gradient of the cost function, the feasibility $\|{X^*}^TAX^*-J\|_F$, the number of iterations, the number of function evaluations, 
					and the elapsed time in second, respectively. In the case that the iteration has not reached the stopping criterion, $X^*$ is replaced by the last iterate. In addition, 
					the convergence histories are displayed in Figure~\ref{fig:tracemin} 
					and 
					Figure~\ref{fig:procrustes} in which we present the value of the cost function, the norm of the Riemannian gradient, and the time elapsed along the iterations, respectively from the left to the right. 
					
					\begin{table}[ht]
						\centering
						\caption{\label{tab:tracemin}Results for minimization of the trace cost function \eqref{eq:trace_fun}.}
					\footnotesize
					\begin{tabular}{l c c c c c c}
						\toprule
						Combination & obj. & grad.& feas. & $\#$iter. & $\#$eval. & CPU \\ 
						\midrule
						Eucl-Cayley & $0.0398$ & $1.0383\e-7$ &$5.1411\e-13$ & 480 & 524 & 750 \\
						
						Eucl-qgeo & $0.0398$ & $1.2258\e-7
						$ &$9.9156\e-08$ & 407 & 446 & 685 \\
						
						gcan-Cayley & 0.0398 & $1.0342\e-7$ & $6.6291\e-13$ & 614 & 666 & 260\\
						
						gcan-qgeo & 0.0398 & $1.1001\e-7$ & $4.8909\e-13$ & 461 & 505 & 324\\
						\bottomrule
					\end{tabular}
				\end{table}
				\begin{table}[ht]
					\centering
					\caption{\label{tab:procrustes_500}Results for minimization of the Procrustes cost function \eqref{eq:procrustes_fun}.} 
				\footnotesize
				\begin{tabular}{l c c c r r r}
					\toprule
					Combination & obj. & grad.& feas. & $\#$iter. & $\#$eval. & CPU \\ 
					\midrule
					Eucl-Cayley & $9.0163\e-4$ & $1.9212\e-3$ &$1.0842\e-12$ & 1000 & 1034 & $1096$ \\
					
					Eucl-qgeo & $7.6519\e-4$ & $2.0001\e-4$ &$8.4243\e-13$ & 1000 & 1034 & $1491$ \\
					
					gcan-Cayley & $2.0092\e-4$ & $2.0370\e-5$ & $4.2354\e-13$ & 435 & 757 & 70\\
					
					
					gcan-qgeo & $1.9463\e-4$ & $1.8871\e-5
					$ & $6.4464\e-13$ & 402 & 734 & 400\\
					\bottomrule
				\end{tabular}
			\end{table}
			
			
			\paragraph{Discussion}
			Through examples presented, we can first see that, all newly established theoretical components work properly and compatibly with the existing ones. Second, as expected, for the cases that $k$ is not small in both models, 
			the schemes using the generalized canonical metric need less time to 
			be convergent with respect to a given \texttt{rstop} compared to the ones using the Euclidean metric although in some cases the required number of iterations can be larger.  Third, the quasi-geodesic retraction apparently is less feasible than the Cayley retraction. This fact might come from the instability of the scaling and squaring algorithm invoked by \texttt{expm} in MATLAB which occasionally 
			occurs as pointed out in \cite{GuetN2016}. 
			And finally, we note that while our sole purpose in constructing the generalized canonical metric is to avoid solving the Lyapunov equation in computing the Riemannian gradient, altering the metric has a great effect on the aspects of numerical methods, such as (numerical) convergence, number of iterations, as one can see in our tests. The result can also depend on the model, namely the cost function. It is advised in  \cite{MishS2016} that the Hessian of the cost function should be taken into account when constructing the metric. To understand the details, we believe that deeper and problem-dependent investigations are required.

			\section{Conclusions}\label{sec:Concl}
			To avoid having to solve the Lyapunov matrix equation during the computation of 
			Riemannian gradient for optimization on the indefinite Stiefel manifold, we proposed 
			the generalized canonical metric and constructed the necessary associated material for running a Riemannian gradient descent scheme, including a retraction based on a quasi-geodesic corresponding to the  new metric. The presented numerical examples verified the theoretical findings and suggested that, using the generalized canonical metric with the Cayley retraction is perhaps the most potential combination since it needs less time to reach a desired approximate solution and delivers 
			better feasibility.
			
			\begin{acknowledgements}
			This research was funded by the TNU-University of Sciences under the research group code: NNC.DHKH.2025.02. 
			\end{acknowledgements}
			
			\subsection*{\textbf{Data availability}} The authors confirm that the data supporting the findings of this study are available within the article.
			\subsection*{\textbf{Competing Interests}} The authors declare that they have no conflict of interest.
				
		\end{document}